\documentclass[11pt]{article}
\usepackage{amssymb,amssymb,amscd}
\begin{document}
\newtheorem{thm}{Theorem}
\newtheorem{lem}[thm]{Lemma}
\newtheorem{cor}[thm]{Corollary}
\newtheorem{prop}[thm]{Proposition}
\def\endpf{\hfill$\bullet$\medskip}
\newtheorem{rem}[thm]{Remark}
\newtheorem{exam}[thm]{Example}
\newtheorem{defn}[thm]{Definition}
\newcounter{jolist}
\newenvironment{joliste}{\begin{list}{{\rm (\roman{jolist})}}{%
    \usecounter{jolist}\setlength{\labelwidth}{10mm}
    \setlength{\leftmargin}{10mm}\setlength{\itemsep}{-4pt}
   \setlength{\topsep}{0pt}}}{\end{list}}
%
%%%%%%%%%%%%%%%%%%%%%%
%%%%%Calligraphic%%%%%%%%%%%%%%%%%
%%%%%%%%%%%%%%%%%%%%%%
\newcommand{\co}{\mathcal{O}}
\newcommand{\ck}{\mathcal{K}}
\newcommand{\cF}{\mathcal{F}}
\newcommand{\cl}{\mathcal{L}}
\newcommand{\ci}{\mathcal{I}}
\newcommand{\cf}{\mathcal{F}}
\newcommand{\cb}{\mathcal{B}}
\newcommand{\bbB}{{\mathbb{B}}}
\newcommand{\cp}{\mathcal{P}}
\newcommand{\ca}{\mathcal{A}}
\newcommand{\cri}{\mathcal{R}}
%%%%%%%%%%%%%%%%%%%%%%
%%%%%operators%%%%%%%%%%%%%%%%%
%%%%%%%%%%%%%%%%%%%%%%
\newcommand{\charc}{{\mathop{\rm Char}\,}}
\newcommand{\Ker}{{\mathop{\rm Ker}}}
\newcommand{\Lie}{{\mathop{\rm Lie}\,}}
\newcommand{\Fr}{{\mathop{\rm Fr}}}
%%%%%%%%%%%%%%%%%%%%%%
%%%%%blackboardbold%%%%%%%%%%%%%%%%%
%%%%%%%%%%%%%%%%%%%%%%
\newcommand{\bp}{{\mathbb P}}
\newcommand{\bn}{{\mathbb N}}
\newcommand{\br}{{\mathbb R}}
\newcommand{\bz}{{\mathbb Z}}
\newcommand{\bq}{{\mathbb Q}}
\newcommand{\bc}{{\mathbb C}}
\newcommand{\bd}{{\mathbb D}}
%%%%%%%%%%%%%%%%%%%%%%
%%%%%GREEK%%%%%%%%%%%%%%%%%
%%%%%%%%%%%%%%%%%%%%%%
\newcommand{\al}{\alpha}
\newcommand{\avee}{\alpha^\vee}
\newcommand{\bvee}{\beta^\vee}
\newcommand{\eps}{\epsilon}
\newcommand{\wt}{\widetilde}
\newcommand{\delp}{\partial^+}
\newcommand{\delm}{\partial^-}
\newcommand{\lam}{{\lambda}}
\newcommand{\Lam}{{\Lambda}}
\newcommand{\usigma}{{\underline{\sigma}}}
\newcommand{\udelta}{{\underline{\delta}}}
\newcommand{\ukappa}{{\underline{\kappa}}}
\newcommand{\ua}{\underline{a}}
\newcommand{\ub}{\underline{b}}
\newcommand{\uc}{\underline{c}}
\newcommand{\ut}{{\underline{\tau}}}
\newcommand{\uk}{{\underline{\kappa}}}
\newcommand{\upi}{{\underline{\pi}}}
\newcommand{\ueta}{{\underline{\eta}}}
\newcommand{\umu}{{\underline{\mu}}}
\newcommand{\unu}{{\underline{\nu}}}
\newcommand{\ulam}{{\underline{\lam}}}
\newcommand{\us}{{\underline{\sigma}}}
\newcommand{\tbar}{{\overline{\tau}}}
\newcommand{\kbar}{{\overline{\kappa}}}
\newcommand{\sbar}{{\overline{\sigma}}}
\newcommand{\ybar}{{\overline{y}}}
\newcommand{\ts}{\tilde{\sigma}}
\newcommand{\tk}{\tilde{\kappa}}
\newcommand{\tta}{{\tilde{\tau}}}
\newcommand{\tw}{\check{w}}
\newcommand{\tS}{\tilde{S}}
%%%%%%%%%%%%%%%%%%%%%%
%%%%%Gothic%%%%%%%%%%%%%%%%%
%%%%%%%%%%%%%%%%%%%%%%
\newcommand{\km}{{\mathfrak{g}}}
\newcommand{\Lb}{{\mathfrak{b}}}
\newcommand{\Lh}{{\mathfrak{h}}}
\newcommand{\Ln}{{\mathfrak{n}}}
\newcommand{\LU}{{\mathfrak{U}}}
\newcommand{\Lu}{{\mathfrak{u}}}
%%%%%%%%%%%%%%%%%%%%%%
%%%%%OVERUNDERline%%%%%%%%%%%%%%%%%
%%%%%%%%%%%%%%%%%%%%%%
\newcommand{\oa}{{\overline{a}}}
\newcommand{\lbar}{{\overline{\ell}}}
\newcommand{\ud}{{\underline{d}}}
\title{Adapted algebras and standard monomials}
\author{P. Caldero$^*$ and P. Littelmann\footnote{This research has been partially supported by the EC TMR network
``Algebraic Lie Representations'', contract no. ERB FMRX-CT97-0100.}}
\maketitle
\begin{abstract} Let $G$ be a complex semisimple Lie group. The aim of this article is to compare two basis for $G$-modules, namely the standard monomial basis and the dual canonical basis. In particular, we give a sufficient condition for a standard monomial to be an element of the dual canonical basis and vice versa.
\end{abstract}
\section*{Introduction}\label{Introduction}

Let $G$ denote a semisimple, simply connected, algebraic group
defined over an algebraically closed field $k$ of arbitrary characteristic.
We fix a Borel subgroup $B$ and a maximal torus $T\subset B$, denote
$W$ the Weyl group of $G$ with respect to $T$. For a dominant
weight $\lam$ let $V(\lam)$ the corresponding Weyl module and let
$Q$ be the parabolic subgroup which normalizes a highest weight vector.
Let $\cl_\lam$ be the corresponding ample line bundle on $G/Q$.
Consider the embedding $G/Q\hookrightarrow \bp(V(\lam))$ given by
the global sections $H^0(G/Q,\cl_\lam)\simeq V(\lam)^*$.

The aim of standard monomial theory is to give a presentation of the ring
$R=\bigoplus_{n\ge 0} H^0(G/Q,\cl_{n\lam})$, which is compatible with
the natural subvarieties of $G/Q$ as for example the Schubert varieties 
$X(\tau)$, the opposite Schubert varieties $X^\kappa$ and the
Richardson varieties $X_\tau^\kappa=X_\tau\cap X^\kappa$, $\tau,\kappa\in W/W_Q$.
More precisely, the aim is to construct a basis $B\subset  H^0(G/Q,\cl_{\lam})$
such that certain ($=$ standard) monomials of degree $n$ in these basis elements form a 
basis $B(n\lam)$ for $H^0(G/Q,\cl_{n\lam})$, the relations provide an algorithm to write a 
non-standard monomial as linear combination of standard monomials, and the
restrictions $\{m\vert_{X_\tau^\kappa}\mid m\in B(n\lam), m\not\equiv 0\hbox{\rm\ on\  }X_\tau^\kappa\}$ of 
the standard monomials to a Richardson variety $X_\tau^\kappa$,
form a basis for $H^0(X_\tau^\kappa,\cl_{n\lam})$. Such a basis has 
turned out to be a powerful tool in the investigation of the geometry 
of Schubert (and related) varieties, see for example \cite{LLM} and \cite{LiS}.

The concept of an adapted algebra has been introduced by the first author in \cite{Cal, Cal1}
in connection with the Berenstein-Zelevinsky conjecture on $q$-commuting elements
of the dual canonical basis $\bc_q[U^-]$. Let $U\subset B$ be a maximal unipotent subgroup,
let $U^-$ be its opposite unipotent subgroup and let $\Ln^-$ be its Lie algebra.
Denote $U_q(\Ln^-)$ the quantized enveloping algebra and let $\cb\subset U_q(\Ln^-)$ be the
canonical basis. As in the classical case, one has a non-degenerate pairing between the quantized 
algebra of regular functions $\bc_q[U^-]$ on $U^-$ and $U_q(\Ln^-)$, so $\bc_q[U^-]$ is naturally
equipped with a canonical basis, the dual $\cb^*\subset \bc_q[U^-]$ of the canonical basis $\cb$.
This basis has many nice properties, for example the natural injection $V_q(\lam)^*\hookrightarrow
\bc_q[U^-]$ is compatible with the basis, i.e., $\cb^*(\lam)=\cb^*\cap V_q(\lam)^*$ is a 
basis for $V_q(\lam)^*$. So the specialization at $q=1$ provides a basis $\cb_1^*(\lam)\subset H^0(G/Q,\cl_{\lam})$,
of which is known that it is compatible with Schubert (and related) varieties.

So on the one hand, $\cb_1^*(\lam)$ looks like a perfect candidate as a starting basis for a standard monomial 
theory. The problem with this approach is that the multiplicative structure of the dual canonical basis
is hardly understood, see \cite{fomzel}. The adapted algebras occur in this context,
they are maximal subalgebras of $\bc_q[U^-]$ which are spanned as a vector space by a subset $S\subset\cb^*$, and 
all elements of $S$ are multiplicative, i.e., with $b,b'\in S$ also $bb'\in S$. In particular, their
product is again an element of the dual canonical basis! The first author has constructed in \cite{Cal} 
for every reduced decomposition $\tw_0$ of the longest word $w_0$ in the Weyl group an adapted algebra
$\ca_{\tw_0}\subset \bc_q[U^-]$.

On the other hand, the construction in \cite{Li$_2$} of a standard monomial theory for the ring $R$
provides in particular a basis for $H^0(G/Q,\cl_{\lam})$, the so-called path vectors $p_\pi$, which 
also has nice representation theoretic features (see for example \cite{LLM,LiS, LCa}), and it is a natural 
question to ask for the relationship between these two bases. Recall that the path model theory 
\cite{Li$_1$, Li$_3$} provides combinatorial model \cite{joseph, kash2} for the crystal base theory, so 
we can index the the elements $b_\pi\in \cb_1^*(\lam)$ and the path vectors $p_\pi$ by L--S paths
of shape $\lam$. 

The main result of this article can be formulated as follows (for a more precise and more
general formulation see section~\ref{pathcanonical}): In the $q$-commutative parts 
of $\cb^*$ given by the adapted algebras $\ca_{\tw_0}$, the two bases coincide up
to normalization. I.e., if $b_\pi\in \cb^*_1(\lam)$ is the specialization of an element of an adapted algebra 
$\ca_{\tw_0}\subset \bc_q[U^-]$ for some  reduced decomposition $\tw_0$ of $w_0$, then $b_\pi=p_\pi$, up 
to normalization (multiplication by a root of unity). Vice versa, if the chain $(\ut)$ of Weyl group elements
in the L--S path $\pi=(\ut,\ua)$ is compatible with some reduced decomposition of $w_0$,
then $b_\pi=p_\pi$, up to normalization. In the last section, we provide some examples to illustrate 
the connection between the path basis and the dual canonical basis.

%%%%%%%%%%%%%%%%%%%%%%%%%%%%%%%%%%%%%%%%%%%%%%%%%%%%%%%%%%%%%%%%%%%%%%%%%%%%%%%%%%%%%%%%%%%%%%%%%%%%%%%%%%%%
%%%%%%%%%%%%%%%%%%%%%%%%%%%%%%%%%%%%%%%%%%%%%%%%%%%%%%%%%%%%%%%%%%%%%%%%%%%%%%%%%%%%%%%%%%%%%%%%%%%%%%%%%%%%
\section{Notation}\label{notation}
%%%%%%%%%%%%%%%%%%%%%%%%%%%%%%%%%%%%%%%%%%%%%%%%%%%%%%%%%%%%%%%%%%%%%%%%%%%%%%%%%%%%%%%%%%%%%%%%%%%%%%%%%%%%
%%%%%%%%%%%%%%%%%%%%%%%%%%%%%%%%%%%%%%%%%%%%%%%%%%%%%%%%%%%%%%%%%%%%%%%%%%%%%%%%%%%%%%%%%%%%%%%%%%%%%%%%%%%%
Let $G$ be a complex semisimple, simply connected algebraic group with Lie algebra $\km$. 
We fix a Cartan decomposition $\km=\Ln\oplus\Lh\oplus\Ln^-$, where $\Lb=\Lh\oplus\Ln$
is a Borel subalgebra of $\km$ with fixed Cartan subalgebra $\Lh$. 
Let $U^-$ be a maximal unipotent subgroup of $G$ such that $\Lie U^-=\Ln$ and let $\bc[U^-]$
be the algebra of regular functions on $U^-$. The left multiplication of $U^-$ on itself 
induces a natural morphism $\phi$ from the enveloping algebra $U(\Ln^-)$ on the algebra of differential
operators on $\bc[U^-]$. This induces a natural non-degenerate pairing on $U(\Ln^-)\times \bc[U^-]$, mapping
$(u,f)$ to the value of $\phi(u)(f)$ at the identity of $U^-$.

This pairing has an analogue in the quantum setting, i.e., one has a natural non-degenerate
pairing on $U_q(\Ln^-)\times \bc_q[U^-]$. So the algebras can be viewed as one being the restricted
dual of the other.

Kashiwara \cite{Ka} and Lusztig \cite{Lu$_1$} have constructed the so-called {\it global crystal} or {\it canonical
basis} $\cb$ of the algebra $U_q(\Ln^-)$. By the {\it dual canonical basis} we mean the basis
$\cb^*$ of $\bc_q[U^-]$ dual to $\cb$ with respect to the pairing above.

Some elements of $\cb^*$ can be easily described, they correspond to so-called extremal weight vectors
in representations. Let $\lam\in P^+$ be a dominant weight and let $V_q(\lam)$ be the corresponding irreducible 
highest weight representation. Fix a highest weight vector $v_\lam$, and for $\tau\in W/W_\lam$ let $v_\tau$ 
be the unique weight vector in $V_q(\lam)$ of weight $\tau(\lam)$ obtained in the following way: 
let $\tau=s_{i_1}\cdots s_{i_r}$ be a reduced decomposition of $\tau$ and set %$v_\tau\in V_q(\lam)$ the
$$
v_\tau=F_{\al_{i_1}}^{(n_1)}\cdots F_{\al_{i_r}}^{(n_r)}v_\lam,
$$
where $n_1=\left<s_{i_2}\cdots s_{i_r}(\lam),\al_{i_1}^\vee\right>,\ldots,n_r=\left<\lam,\al_{i_r}^\vee\right>$. 
It follows from the Verma relations for the generators of $U_q(\Ln^-)$ that the vector is
independent of the chosen decomposition. This weight space is one dimensional, let $b_\tau^\lam\in V_q(\lam)^*$ be the unique
element of weight $-\tau(\lam)$ such that 
\begin{eqnarray}
b^\lambda_\tau(v_\tau)=1.
\end{eqnarray}
This linear form is considered as a form  on $U_q(\Ln^-)$ by setting $b_\tau^\lam(u)=b_\tau^\lam(u v_\lam)$. 
The compatibility property of the canonical basis with highest weight representations shows 
immediately that these functions $b_\tau^\lam$ are elements of $\cb^*$. In the following we refer
to the $b_\tau^\lam$ as {\it $(\lam,q)$--minors}. For $\km={\mathfrak sl}_n$
and $\lam$ a fundamental weight, these are exactly the quantum minors.

It is often more convenient to forget in the notation the stabilizer $W_\lam$ and, by abuse of notation, just
write $b_\tau^\lam$ for $\tau\in W$ and $\lam\in P^+$. In this way we can formulate the
following simple product rule (see \cite{DP}):
\begin{eqnarray}\label{minorproducts}
\forall\lam,\mu\in P^+:\quad b_\tau^\lam b_\tau^\mu =b_\tau^{\lam+\mu}\quad\hbox{\rm\ in\ }\bc_q[U^-].
\end{eqnarray}
More generally, for $\xi\in V_q(\lam)^*$ let $c_\xi^\lam$ be the linear form on $U_q(\Ln^-)$ obtained by setting
$c_\xi^\lam(u):=c_\xi^\lam(u v_\lam)$. The following commutation relation can be found for example in \cite{DP},
it is a direct consequence of \cite{LeS}, where the expression of the ${\cal R}$-matrix is made explicit:
suppose $\xi,\eta$ are of weight $\nu_1,\nu_2$. Then
\begin{eqnarray}\label{generalproducts}
c^\lam_{\xi} c^\lam_{\eta}= q^{-(\lam,\lam)+(\nu_1,\nu_2)}c^\lam_{\eta}c^\lam_{\xi} 
+\sum c^\lam_{\eta_i}c^\lam_{\xi_i}
\end{eqnarray}
where there exist $p_i\in\bc(q)$ and some non-constant monomials $M_i$ such that
$$
\eta_i\otimes\xi_i=p_i(q)(M_i(E)\otimes M_i(F)) \eta\otimes\xi.
$$
%%%%%%%%%%%%%%%%%%%%%%%%%%%%%%%%%%%%%%%%%%%%%%%%%%%%%%%%%%%%%%%%%%%%%%%%%%%%%%%%%%%%%%%%%%%%%%%%%%%%%%%%%%%%
%%%%%%%%%%%%%%%%%%%%%%%%%%%%%%%%%%%%%%%%%%%%%%%%%%%%%%%%%%%%%%%%%%%%%%%%%%%%%%%%%%%%%%%%%%%%%%%%%%%%%%%%%%%%
\section{$(\lam,q)$--minors and adapted algebras}\label{qminors}
%%%%%%%%%%%%%%%%%%%%%%%%%%%%%%%%%%%%%%%%%%%%%%%%%%%%%%%%%%%%%%%%%%%%%%%%%%%%%%%%%%%%%%%%%%%%%%%%%%%%%%%%%%%%
%%%%%%%%%%%%%%%%%%%%%%%%%%%%%%%%%%%%%%%%%%%%%%%%%%%%%%%%%%%%%%%%%%%%%%%%%%%%%%%%%%%%%%%%%%%%%%%%%%%%%%%%%%%%
Recall that two elements $a,b\in \bc_q[U^-]$
are said to {\it $q$-commute} if $ab=q^m ba$ for some $m\in\bz$. Two elements $b,b'\in\cb^*$ are called 
{\it multiplicative} if the product is a multiple of an element of the canonical basis, i.e., $bb'=rb''$ 
for some $b''\in\cb^*$, $r\in\bc(q)$. The elements of the dual canonical basis are in general neither 
multiplicative nor $q$-commuting. 
Let $\varpi_1,\ldots,\varpi_n$ be the set of fundamental weights and fix a reduced decomposition 
$w_0=s_{i_1}\cdots s_{i_N}$ of the longest word in $W$, we write $\tw_0$ to refer to this decomposition. Set 
\begin{eqnarray}\label{generatorforadap}
y_j=s_{i_1}\cdots s_{i_j} \hbox{\rm\ and set\ }
%b_j\in\cb^*\hbox{\rm\  be the element\ }
b_j:= b_{y_j}^{\varpi_{i_j}}\in\cb^*.
\end{eqnarray}

\begin{rem}\label{modulodescription}\rm
Since $b_\tau^\lam=b_{\tau'}^\lam$ for $\tau\equiv\tau'\bmod W_\lam$, it is easy to verify the following
description of the elements above:
\[
\{b_1,\ldots,b_N\}=\{ b_{y_j}^{\varpi_i}\mid 1\le i\le n,1\le j\le N\}
\]
\end{rem}

These elements have been studied in great detail by the first author in \cite{Cal}, where he shows that 
they have the following remarkable properties:
\begin{joliste}
\item[{\rm i)}] they are multiplicative and $q$-commute, and, more precisely, 
\item[{\rm ii)}] the monomials $b_1^{m_1}\cdots b_N^{m_N}$ are, up to a power of $q$, elements of $\cb^*$.
\end{joliste}
Let $S_{\tw_0}$ be the set of the monomials above and let $\ca_{\tw_0}$ be the subalgebra of $\bc_q[U^-]$
spanned by $S_{\tw_0}$. Actually, ${\cal A}_{\tilde w_0}$ is generated as a space by $S_{\tilde w_0}$.

These properties motivate the following definition: a subalgebra $\ca$ of $\bc_q[U^-]$ is called {\it adapted}
if it is spanned by a subset $\cp^*\subset\cb^*$ such that the elements in $\cp^*$ are multiplicative
and $\ca$ is maximal with this property, i.e., for all $b\in\cb^*-\cp^*$ there exist a $p\in \cp^*$ such that
$b$ and $p$ are not multiplicative.

\begin{thm}[\cite{Cal}]\label{adaptedalgebrathm} The subalgebra $\ca_{\tw_0}$ is adapted with the elements of $S_{\tw_0}$
(up to a rescaling by a power of $q$) as spanning set. The set $S_{\tw_0}$ is an Ore set
in $\bc_q[U^-]$, and for the localizations one has: $S^{-1}_{\tw_0}\bc_q[U^-]=S^{-1}_{\tw_0}\ca_{\tw_0}$. 
\end{thm}
%%%%%%%%%%%%%%%%%%%%%%%%%%%%%%%%%%%%%%%%%%%%%%%%%%%%%%%%%%%%%%%%%%%%%%%%%%%%%%%%%%%%%%%%%%%%%%%%%%%%%%%%%%%%
%%%%%%%%%%%%%%%%%%%%%%%%%%%%%%%%%%%%%%%%%%%%%%%%%%%%%%%%%%%%%%%%%%%%%%%%%%%%%%%%%%%%%%%%%%%%%%%%%%%%%%%%%%%%
\section{Lusztig's Frobenius maps}\label{frobeniusmaps}
%%%%%%%%%%%%%%%%%%%%%%%%%%%%%%%%%%%%%%%%%%%%%%%%%%%%%%%%%%%%%%%%%%%%%%%%%%%%%%%%%%%%%%%%%%%%%%%%%%%%%%%%%%%%
%%%%%%%%%%%%%%%%%%%%%%%%%%%%%%%%%%%%%%%%%%%%%%%%%%%%%%%%%%%%%%%%%%%%%%%%%%%%%%%%%%%%%%%%%%%%%%%%%%%%%%%%%%%%

Let $\cri$ be the ring of Laurent polynomials $\bz[q,q^{-1}]$. We denote 
$U_\cri(\Ln^-)$ the $\cri$-form of $U_q(\Ln^-)$ generated by the 
divided powers $F_i^{(n)}$. Fix a primitive $2\ell$-th root of unity $\xi$,
where $\ell$ itself is even, and denote $\bz_\xi$ the ring $\bz[\xi]$.

Let $U(\Ln^-)$ be the $\bz_\xi$-form of the enveloping algebra of $\Ln^-$,
i.e., it is obtained from the Kostant $\bz$-form by the ring extension $\bz\hookrightarrow
\bz_\xi$. 

Let $U_\xi(\Ln^-)$ be the $\bz_\xi$-algebra obtained from $U_\cri(\Ln^-)$ by specialization 
at $q=\xi$, i.e., it is obtained from the Lusztig $\cri$-form $U_\cri(\Ln^-)$ by base
change with respect to $\cri\rightarrow \bz_\xi$, $q\mapsto \xi$.

Let $A=(a_{i,j})_{1\le i,j\le n}$ be the Cartan matrix of $G$ and denote by 
$^tA=(\oa_{i,j})$, $\oa_{i,j}:=a_{j,i}$, the transposed matrix. Let $\ud=(d_1,\ldots,d_n)$, $d_i\in\bn$,  
be minimal such that the matrix $(d_i a_{i,j})$ is symmetric. We denote by $d$ the smallest common multiple 
of the $d_j$. To distinguish between the objects associated to $A$ and $^tA$, we add always
a $^t$ in the notation. For example, we write ${}^t\km={}^t\Ln^-\oplus{}^t\Lh\oplus{}^t\Ln$ for the
Cartan decomposition of the semisimple Lie algebra associated to ${}^tA$, and we write
${}^tF_i$ for the generators of $U_q({}^t\Ln^-)$.

In the following we assume that $\ell$ is divisible by $2d$ and we set $\lbar=\ell/d$
and let $\ell_i\in\bn$ be minimal such that $\ell_i \frac{d}{d_i}\in\ell\bz$.
In \cite{Lu$_1$}, Chapter 35, Lusztig defines
two algebra homomorphisms:
\[
\begin{array}{rclccccccc}
\Fr   &:& {}^tU_\xi(\Ln^-)&\rightarrow& U(\Ln^-)       &\quad &{}^tF_i^{(m)}&\mapsto& F_i^{(m/\ell_i)}\\
{\Fr'}&:& U(\Ln^-)        &\rightarrow&{}^tU_\xi(\Ln^-)&\quad &    F_i^{(m)}&\mapsto& {}^tF_i^{(\ell_i m)}
\end{array}
\]
where we set $F_i^{(m/\ell_i)}:=0$ if $\ell_i$ does not divide $m$. The
first map is in fact a Hopf algebra homomorphism, but the second is not.

The corresponding $\cri$-form of $\bc_q[U^-]$ is denoted $\cri_q[U^-]$, and its specializations 
at $\xi$ respectively $1$ are denoted $\cri_\xi[U^-]$ respectively $\cri_1[U^-]$. 
The two algebra homomorphisms $\Fr$ and $\Fr'$ induce dual maps:
\begin{eqnarray*}\label{dualfrobenius}
\Fr^*: \cri_1[U^-]\rightarrow  \cri_\xi[{}^tU^-]&\quad&{\Fr'}^*:  \cri_\xi[{}^tU^-]\mapsto  \cri_1[U^-].
\end{eqnarray*}
Note that $\Fr^*$ is an algebra homomorphism since $\Fr$ is a Hopf algebra homomorphism.
Further, since $\Fr\circ\Fr'$ is the identity map it follows that ${\Fr'}^*\circ\Fr^*$ is the identity.
In particular, $\Fr^*$ is injective and maps $\cri_1[U^-]$ isomorphically onto a commutative subalgebra of
$\cri_\xi[{}^tU^-]$

The maps $\Fr$ respectively its dual $\Fr^*$ are called the {\it Frobenius maps}, and the maps
${\Fr'}^*$ and ${\Fr'}^*$ are called the {\it Frobenius splittings} (for the connection
with the algebraic geometric notion of Frobenius splitting in positive characteristic
see \cite{KL1} and \cite{KL2}).

Of special interest are for us the $(\lam,q)$--minors:
\begin{prop}\label{frobeniusminors}
The Frobenius map $\Fr^*$ is an algebra homomorphism which, restricted to the $(\lam,q)$--minors, is 
the $\lbar$-th power map, i.e., 
$$
\Fr^*(b_\tau^\lam)=(b_\tau^{\lam})^\lbar=b_\tau^{\lbar\lam}.
$$ 
The splitting ${\Fr'}^*$ is a $\bz_\xi$-module homomorphism which commutes with the image
of $\Fr^*$, i.e., for $f\in \cri_1[U^-]$ and $g\in \cri_\xi[{}^tU^-]$ we have
\[
{\Fr'}^*( {\Fr}^*(f)g)=f {\Fr'}^*(g).
\]
In particular,
$$
{{\Fr}'}^*(b_{\tau_1}^{m_1\lbar}\cdots b_{\tau_r}^{m_r\lbar})=b_{\tau_1}^{m_1}\cdots b_{\tau_r}^{m_r}.
$$
\end{prop}

\par\noindent
{\it Proof.} The first part follows immediately from \cite{Li$_2$}, we sketch for convenience the
arguments: let $V_\xi(\lbar\lam)$ be the $U_{\xi}({}^t\km)$--Weyl module for the highest weight
$\lbar\lam$, and let $L(\lbar\lam)$ be the simple module of highest weight $\lbar\lam$. Then
$L(\lbar\lam)$ is nothing but the Weyl module $V(\lam)$ for $U(\km)$, viewed as quantum module 
via the (extension of the) Frobenius map $\Fr:U_{\xi}({}^t\km)\rightarrow U(\km)$. One sees easily 
that this induces a dual map (also) denoted $\Fr^*: V(\lam)^*\rightarrow V_\xi(\lbar\lam)^*$
which sends extremal weight vectors of $V(\lam)^*$ to the corresponding extremal weight vectors 
of $V_\xi(\lbar\lam)^*$. In particular, the $b_\tau^\lam$ is sent to $b_\tau^{\lbar\lam}$, which is 
equal to $(b_\tau^{\lam})^\lbar$ by (\ref{minorproducts}). 
The commutation part follows from \cite{KL1}, see also \cite{KL2}. More precisely, the assumptions made on $\ell$ 
in these articles are, as mentioned in their introductions, just made for convenience. Suppose $f\in \cri_1[U^-]$ and 
$g\in \cri_\xi[{}^tU^-]$. By choosing $\lam\in P^+$ generic enough, we can assume that 
we can identify $f$ and $g$ as elements of representations, i.e.,
\[
f\in V(\lam)^*\hookrightarrow \cri_1[U^-] ,\quad g\in V_\xi(\lbar\lam)^*\hookrightarrow \cri_\xi[{}^tU^-],
\] 
and ${\Fr'}^*(\Fr^*(f)g)\in  V(2\lam)^*\hookrightarrow \cri_1[U^-]$. Now the same
arguments as in \cite{KL1}, section~4,  or \cite{KL2}, section~6, go through to prove the compatibility
${\Fr'}^*(\Fr^*(f)g)=f{\Fr'}^*(g)$. The rest is obvious since ${\Fr'}^*$ is a splitting.
\endpf

The canonical basis $\cb\subset U_q(\Ln^-)$ spans $U_\cri(\Ln^-)$ as an $\cri$-module, correspondingly
its dual basis $\cb^*$ spans $\cri_q[U^-]$ as an $\cri$-module. It follows that the images of this basis
in $\cri_\xi[{}^tU^-]$ respectively $\cri_1[U^-]$ span the corresponding specializations.

Fix a reduced decomposition $w_0=s_{i_1}\cdots s_{i_N}$ of the longest word in $W$. As in section~\ref{qminors}, 
set $b_j:= b_{y_j}^{\varpi_{i_j}}\in\cb^*\subset U_q(\km)$, let $S_{\tw_0}$ be the set of the monomials in the $b_j$. 
Still denote  ${\cal A}_{\tilde w_0}$ the  ${\cal
R}$-module generated by $S_{\tilde w_0}$.
%and let $\ca_{\tw_0}$ be the subalgebra of $\bc_q[U^-]$ spanned by $S_{\tw_0}$. 
The ${\cal R}$--algebra $\ca_{\tw_0}$ specializes
to a subalgebra $\ca^\xi_{\tw_0}\subset\cri_\xi[{}^tU^-]$ respectively $\ca^1_{\tw_0}\subset\cri_1[U^-]$.
In both cases the algebra has the (images of the) monomials in $S_{\tw_0}$ as a basis as $\bz_\xi$-module.

As an immediate consequence of the proposition above one sees:
\begin{cor}\label{frobeniusadapted} 
$\Fr^*$ induces an inclusion of algebras $\Fr^*:\ca^1_{\tw_0}\rightarrow {}^t\ca^\xi_{\tw_0}$ 
such that $\Fr^*(b_1^{m_1}\cdots b_N^{m_N})=b_1^{m_1\lbar}\cdots b_N^{m_N\lbar}$. 
%The splitting is a $\bz_\xi$-module homomorphism ${\Fr'}^*:\cri_\xi[U^-]\rightarrow \cri_1[U^-]$ such that
%${\Fr'}^*(b_1^{m_1\lbar}\cdots b_N^{m_N\lbar})=b_1^{m_1}\cdots b_N^{m_N}$.
\end{cor}

%%%%%%%%%%%%%%%%%%%%%%%%%%%%%%%%%%%%%%%%%%%%%%%%%%%%%%%%%%%%%%%%%%%%%%%%%%%%%%%%%%%%%%%%%%%%%%%%%%%%%%%%%%%%
%%%%%%%%%%%%%%%%%%%%%%%%%%%%%%%%%%%%%%%%%%%%%%%%%%%%%%%%%%%%%%%%%%%%%%%%%%%%%%%%%%%%%%%%%%%%%%%%%%%%%%%%%%%%
\section{The path model, $(\lam,q)$--minors and SMT}\label{pathmodel}
%%%%%%%%%%%%%%%%%%%%%%%%%%%%%%%%%%%%%%%%%%%%%%%%%%%%%%%%%%%%%%%%%%%%%%%%%%%%%%%%%%%%%%%%%%%%%%%%%%%%%%%%%%%%
%%%%%%%%%%%%%%%%%%%%%%%%%%%%%%%%%%%%%%%%%%%%%%%%%%%%%%%%%%%%%%%%%%%%%%%%%%%%%%%%%%%%%%%%%%%%%%%%%%%%%%%%%%%%
We come now first back to representations of $G$ respectively $U(\km)$. A combinatorial
tool for the analysis of these representations is the path model, of which
we recall quickly the most important features. Let $\lam\in P^+$ be a dominant weight. 

Let $\ut=(\tau_0,\ldots,\tau_r)$ be a strictly decreasing sequence (with respect to the 
Bruhat order) of elements of $W/W_\lam$, and let $\ua=(a_1,\ldots,a_r)$ 
be a strictly increasing sequence of rational numbers such that $0<a_1<\ldots<a_r<1$. 
The pair $\pi=(\ut,\ua)$ is called a {\it L-S path of shape $\lam$} if the pair satisfies the 
following integrality condition. For all $i=1,\ldots,r$:
\begin{itemize}
\item set ${r_i}=l(\tau_{i-1})-l(\tau_i)$. There exists a sequence $\beta_1,\ldots,\beta_{r_i}$ 
of positive roots joining $\tau_{i-1}$ and $\tau_i$ by the corresponding reflections, i.e.,
\[
\tau_{i-1} >s_{\beta_1}\tau_{i-1}>s_{\beta_2}s_{\beta_1}\tau_{i-1}>\ldots
>s_{\beta_{r_i}}\cdots s_{\beta_1}\tau_{i-1}=\tau_i,
\]
and $a_i\left<\tau_{i-1}(\lam),s_{\beta_1}\cdots s_{\beta_{j-1}}(\beta^\vee_j)\right>\in\bz$ for all
$j=1,\ldots,r_i$.
\end{itemize}
In the following formulas we set always $a_0=0$ and $a_{r+1}=1$.
The {\it weight of an L-S path} $\pi=(\ut,\ua)$ of shape $\lam$ is the convex linear combination
$$
\pi(1)=\sum_{i=0}^r x_{i+1}\tau_{i}(\lam), \hbox{\rm\ where\ }x_i=a_i-a_{i-1}\hbox{\rm\ for\ } 1\le i\le r+1.
$$
For more details on the combinatorics of L-S paths we refer to \cite{Li$_1$}. 
Let $B(\lam)$ be the set of L-S paths of shape $\lam$, and let $B(\lam)_\tau$
be the set of L-S paths $\pi=(\ut,\ua)$ of shape $\lam$ such that $\tau\ge\tau_0$.
The character of the Weyl module $V(\lam)$ of highest weight 
$\lam$ and the Demazure module $V(\lam)_\tau$ can be calculated using the L-S paths: 

\begin{thm}[\cite{Li$_1$}]\label{pathcharacter} 
\[
\charc V(\lam)=\sum_{\pi\in B(\lam)} e^{\pi(1)}\hbox{\rm\ and\ }
\charc V(\lam)_\tau=\sum_{\pi\in B(\lam)_\tau} e^{\pi(1)}.
\]
\end{thm}

The path model of a representation is closely connected to the combinatorics of the 
{\it crystal base}. For the irreducible representation $V_q(\lam)$ fix a highest
weight vector $v_\lam$ and let $\cb(\lam)=\{bv_\lam\mid b\in\cb, b v_\lam\not=0\}$
be the canonical basis of $V_q(\lam)$. Denote $L_q(\lam)$ the $\bc[q]$--submodule
of $V_q(\lam)$ spanned by the elements of the canonical basis $\cb(\lam)$ and
set 
$$
B^q(\lam)=\{ b \bmod L_q(\lam)\mid b\in \cb(\lam)\}.
$$
Recall (\cite{kash2,joseph}) that the Kashiwara crystal $B^q(\lam)$ (with the Kashiwara 
operators $f_\al,e_\al$) is isomorphic to the crystal given by the path model $B(\lam)$ 
(with the root operators). In the following we will often identify these two 
crystals and just write $B(\lam)$. For an element $b\in \cb(\lam)$ we
write $[b]^\lam$ or just $[b]$ for its class in $B(\lam)$, and for an element
$\pi\in B(\lam)$ we write $b_\pi$ for the element in $\cb(\lam)$ such that $[b_\pi]=\pi$.

The tensor products of crystals can be defined in the path model in terms of
concatenations of paths or in terms of {\it standard tuples}:

\begin{defn}\label{definingchain}\rm
Let $\pi_1=(\ut,\ua)\in B(\lam_1)$, $\pi_2=(\uk,\ub)\in B(\lam_2),\ldots$, respectively 
$\pi_m=(\us,\uc)\in B(\lam_m)$ be L-S paths of shape $\lam_1$,
$\lam_2,\ldots$, respectively $\lam_m$, where
\[
\ut=(\tau_0,\ldots,\tau_r),\quad \uk=(\kappa_0,\ldots,\kappa_s),\ldots,\quad\us=(\sigma_0,\ldots,\sigma_t).
\]
The tuple $\upi=(\pi_1,\pi_2,\ldots,\pi_m)$ is called a {\it standard tuple of shape} $\ulam=(\lam_1,\ldots,\lam_m)$ 
if there exist representatives $\tta_0,\ldots\tta_r\in W$ of the $\tau_0,\ldots,\tau_r\in W/W_{\lam_1}$,
and representatives $\tk_0,\ldots,\tk_s\in W$ of the $\kappa_0,\ldots,\kappa_s\in W/W_{\lam_2}$ etc., such that
\[
\tta_0\ge\ldots\ge\tta_r\ge\tk_0\ge\ldots\ge\tk_s\ge\ldots\ge\ts_0\ge\ldots\ge\ts_t,
\]
Such a sequence $(\tta_0,\ldots,\tta_r,\tk_0,\ldots,\tk_s,\ldots,\ts_0,\ldots,\ts_t)$ of representatives 
is called a {\it defining chain}. The weight $\upi(1)$ is the sum $\pi_1(1)+\ldots+\pi_m(1)$ of the weights
of the L--S paths.
\end{defn}

\begin{rem}\rm
If $\lam_1=\ldots=\lam_m$, the condition on standardness reduces to the condition
$\tau_r\ge\kappa_0\ge\kappa_s\ge\ldots\ge\sigma_0$ in the Bruhat ordering on $W/W_\lam$.
\end{rem}

By the independence of the crystal structure of path models \cite{Li$_3$}, we have a canonical bijection
between the set $B(\ulam)$ of standard tuples of shape $\ulam$ and the set of
L--S paths $B(\lam)$ of shape $\lam=\lam_1+\ldots+\lam_m$. In fact, the standard tuples
correspond exactly to the elements in the Cartan component of the tensor product of the crystals.
The notation $b_\upi\in \cb(\lam)$ for an element of the canonical basis corresponding to an element of
a path model generalizes in the obvious way.

The following is a translation of \cite{LNT}, Proposition 33, into the language of standard tuples:

\begin{prop}\label{lnt}
Let $b_i\in\cb(\lam_i)$ be elements of the canonical basis and let $\pi_i\in B(\lam_i)$
be the L--S path such that $[b_i]=\pi_i$, $i=1,\ldots,m$. Suppose the tuple $\upi=(\pi_1,\ldots,\pi_m)$ is
standard, let $\pi\in B(\lam)$, $\lam=\sum_{i=1}^m\lam_i$, be the L--S path of shape
$\lam$ that identifies with the standard tuple, and let $b_\pi=b_\upi\in\cb(\lam)$ be such that $[b_\pi]=\pi$.
Then $[b_1]\otimes\ldots\otimes[b_m]$ identifies with $[b_\pi]$, and for the dual basis elements
there exists a $s\in\bn$ such that 
$$
q^s b_1^*b_2^*\ldots b_m^*=b_\pi^* \bmod qL_q^*(\lam).
$$
\end{prop}

The dual representation $V(\lam)^*$ can be geometrically realized as the space of global 
sections $H^0(G/B,\cl_\lam)$ of the line bundle $\cl_\lam$ on $G/B$.
We recall now the construction of a basis of $H^0(G/B,\cl_\lam)$ using $(\lam,q)$--minors.
The restriction map induces a map between (tensor products of) dual Weyl modules of the 
quantum group $U_\xi({}^t\km)$:
\[
V_\xi(\lam)^*\otimes\ldots\otimes V_\xi(\lam)^*\longrightarrow V_\xi(\lbar\lam)^*.
\]
We write shortly $b_\tau^\lam\cdots b_\kappa^\lam$ for the image of 
$b_\tau^\lam\otimes\ldots\otimes b_\kappa^\lam$. The dual of the Frobenius splitting 
induces a map between dual Weyl modules of the quantum group $U_\xi({}^t\km)$
and $U(\km)$ (see \cite{Li$_2$}):
\[
{\Fr'}^*:V_\xi(\lbar\lam)^* \longrightarrow V(\lam)^*.
\]
Note that these maps are actually defined over $\bz_\xi$.

Let now $\pi=(\ut,\ua)\in B(\lam)$ be an L-S path of shape $\lam$, fix $\ell\in \bn$, $\ell>0$ minimal
such that $2d$ divides $\ell$ and $\lbar a_i\in\bn$ for all $i=1,\ldots,r$, and
let $\xi$ be a corresponding primitive $2\ell$-th root of unity. 

\begin{defn}\label{pathvecdefin}
The path vector $p_\pi\in H^0(G/B,\cl_\lam)$, $\pi=(\ut,\ua)$, is defined as:
\[
p_\pi:={\Fr'}^*\left( (b_{\tau_0}^\lam)^{\lbar x_1}(b_{\tau_{1}}^\lam)^{\lbar x_2}\cdots 
(b_{\tau_{r-1}}^\lam)^{\lbar x_{r}}(b_{\tau_{r}}^\lam)^{\lbar x_{r+1}}\right),
\]
where $x_i=a_i-a_{i-1}$ for $1\le i\le r+1$, and the $b_{\tau_j}^\lam$ are $(\lam,q)$-minors in $V_\xi(\lam)^*$ for $j=0,\ldots,r$.
\end{defn}

To distinguish between extremal weight vectors for $U(\km)$ and $U_\xi({}^t\km)$ we keep the notation
$b_\tau^\lam$ in the quantum case and write $p_\tau^\lam$ for the classical case.
The reader should think of the $p_\pi$ as a kind of algebraic
approximation of
$$
p_\pi\sim {}^{\ell}\sqrt{p_{\tau_0}^{\ell x_1}p_{\tau_1}^{\ell x_2}p_{\tau_3}^{\ell x_3}\ldots
p_{\tau_r}^{\ell x_{r+1}}}.
$$
Though the expression above does not make sense in general, we will see soon that at least in some 
special cases the expression can be given a useful interpretation. 

We will also consider {\it standard monomials} 
in the path vectors:

\begin{defn}\rm
Let $\upi=(\pi_1,\ldots,\pi_m)$ be a tuple of shape $\ulam=(\lam_1,\ldots,\lam_m)$. We associate to
$\upi$ the monomial $p_\upi=p_{\pi_1} p_{\pi_2}\cdots p_{\pi_m}$ in $H^0(G/Q,\cl_{\lam})$,
where $\lam=\lam_1+\ldots+\lam_m$. The monomial $p_\upi$ is called a {\it standard monomial of shape} $\ulam$
if $\upi$ is a standard tuple.
\end{defn}

The path vectors are defined over $\bz_\xi$ for some appropriate root of unity,
so the collection ${\bbB}(\ulam)$ of all standard monomials of 
shape $\ulam$ is defined over the ring $S\subset\bc$ generated by $\bz$ and all roots of unity.

\begin{thm}[\cite{Li$_2$}]
The set  of standard monomials $\bbB(\ulam)$ of shape $\ulam$ forms a basis
for the space of sections $H^0(G/B,\cl_{\lam})$, $\lam=\lam_1+\ldots+\lam_m$.
\end{thm}

For more information about applications of standard monomial theory (singularities of Schubert varieties,
generators and relations for homogeneous coordinate rings of Schubert varieties, Koszul property, 
Pieri-Chevalley type formula, ...) see for example \cite{LLM}, \cite{Li$_2$}, \cite{LiS}.

%%%%%%%%%%%%%%%%%%%%%%%%%%%%%%%%%%%%%%%%%%%%%%%%%%%%%%%%%%%%%%%%%%%%%%%%%%%%%%%%%%%%%%%%%%%%%%%%%%%%%%%%%%%%
%%%%%%%%%%%%%%%%%%%%%%%%%%%%%%%%%%%%%%%%%%%%%%%%%%%%%%%%%%%%%%%%%%%%%%%%%%%%%%%%%%%%%%%%%%%%%%%%%%%%%%%%%%%%
\section{The path vectors and the dual canonical basis}\label{pathcanonical}
%%%%%%%%%%%%%%%%%%%%%%%%%%%%%%%%%%%%%%%%%%%%%%%%%%%%%%%%%%%%%%%%%%%%%%%%%%%%%%%%%%%%%%%%%%%%%%%%%%%%%%%%%%%%
%%%%%%%%%%%%%%%%%%%%%%%%%%%%%%%%%%%%%%%%%%%%%%%%%%%%%%%%%%%%%%%%%%%%%%%%%%%%%%%%%%%%%%%%%%%%%%%%%%%%%%%%%%%%
Given an L-S path of shape $\lam$, we use the same notation $p_\pi$ for the path vector in $H^0(G/B,\cl_\lam)$
as well as the element $c_{p_\pi}^\lam\in \bc[U^-]$. 
The dual canonical basis $\cb^*$ is compatible with the injection $H^0(G/B,\cl_\lam)\hookrightarrow \bc[U^-]$,
i.e., a subset of $\cb^*$ forms a basis of the image of $H^0(G/B,\cl_\lam)$. 

So it is natural to ask for a description of the transformation matrix between the basis given by the path vectors
and the basis given by elements of the dual canonical basis.

Let $\pi=(\ut,\ua)$ be an L-S path of shape $\lam$. We identify an element $\tau\in W/W_\lam$
with its minimal representative in $W$. Fix a reduced decomposition $w_0=s_{i_1}\cdots s_{i_N}$ 
of the longest word of the Weyl group, we write $\tw_0$ to refer to this decomposition. 
For $j=1,\ldots, N$, set $y_j=s_{i_1}\cdots s_{i_j}$ as in (\ref{generatorforadap}).

\begin{defn}\rm
We say that {\it $\pi$ is compatible with $\tw_0$} if $\ut\subset \{y_1,\ldots,y_N\}$, i.e.,
$\{\tau_0,\ldots,\tau_r\}\subset\{y_1,\ldots,y_N\}$.
\end{defn}

Let now $\upi=(\pi_1,\ldots,\pi_m)$ be a standard tuple of shape $\ulam=(\lam_1,\ldots,\lam_m)$ with a defining chain
(see Definition~\ref{definingchain}).
\[
\tta_0\ge\ldots\ge\tta_r\ge\tk_0\ge\ldots\ge\tk_s\ge\ldots\ge\ts_0\ge\ldots\ge\ts_t,
\]
\begin{defn}\rm
We say that {\it $\upi$ is compatible with $\tw_0$} if the defining chain can be chosen
such that  $\{\tta_0,\ldots,\tta_r,\tk_0,\ldots,\tk_s,\ldots,\ts_0,\ldots,\ts_t\}\subset\{y_1,\ldots,y_N\}$.
\end{defn}

The following theorem provides the connection between certain standard monomials and the adapted algebras:

\begin{thm}\label{pathandcanonicalzigzag}
Let $p_\upi=p_{\pi_1}\cdots p_{\pi_m}\in H^0(G/B,\cl_{\lam})$ be a standard monomial monomial
of shape $\ulam=(\lam_1,\ldots,\lam_m)$. If there exists 
a defining chain $(\tta_0,\ldots,\tta_r,\ldots,\ts_0,\ldots,\ts_t)$ which is compatible with some reduced 
decomposition of $w_0$, then the inclusion  $H^0(G/B,\cl_{\lam})\hookrightarrow \bc[U^-]$,
$\lam=\lam_1+\ldots+\lam_m$, maps the standard 
monomial $p_\upi$ onto the element $b^*_\upi$ of the dual canonical basis, up to multiplication by a root of unity.
\end{thm}

\noindent
{\it Proof.} We fix a reduced decomposition $\tw_0$ of $w_0$ and assume that the
standard tuple $\upi=(\pi_1,\ldots,\pi_m)$ has a defining chain 
compatible with $\tw_0$. Let $\ca_{\tw_0}^1\subset\cri_1[U^-]\subset \bc[U^-]$
be the corresponding (specialization of the) adapted algebra defined in section~\ref{qminors}. 
We will show that the $p_{\pi_j}$ are elements of the set of monomials $S_{\tw_0}$ spanning 
$\ca_{\tw_0}^1$. This proves the theorem by Theorem~\ref{adaptedalgebrathm}.

So in the following it suffices to consider only one path vector $p_\pi$. Fix an appropriate
$\ell\in \bn$ as in section~\ref{pathmodel}, Definition~\ref{pathvecdefin}, for
$\pi=(\tau_0,\ldots,\tau_r;a_1,\ldots,a_r)$, and let $\xi$ be a primitive $2\ell$-th root of
unity. Set
\[
m_\pi:=(b_{\tau_0}^\lam)^{\lbar a_1}(b_{\tau_{1}}^\lam)^{\lbar (a_2-a_1)}\cdots 
(b_{\tau_{r-1}}^\lam)^{\lbar(a_r-a_{r-1})}(b_{\tau_{r}}^\lam)^{\lbar(1-a_r)},
\]
where the $b_{\tau_j}^\lam$ are $(\lam,q)$--minors in $V_\xi(\lam)^*$ for $j=0,\ldots,r$.
By definition, we have $p_\pi={\Fr'}^*(m_\pi)$. Let $\lam=\sum_{j=1}^n\lam_j\varpi_j$ be the
expression of $\lam$ as a linear combination of fundamental weights. Then, by (\ref{minorproducts}), 
we know for arbitrary $\tau\in W$:
$$
b_{\tau}^\lam=\prod_{i=1}^n (b_{\tau}^{\varpi_j})^{\lam_j}
$$
So we can rewrite $m_\pi$ as:
\[
m_\pi=\left(\prod_{j=1}^n (b_{\tau_0}^{\varpi_j})^{\lbar a_1\lam_j}\right)
\left(\prod_{i=j}^n (b_{\tau_1}^{\varpi_j})^{\lbar (a_2-a_1)\lam_j}\right)
\cdots
\left(\prod_{i=j}^n (b_{\tau_r}^{\varpi_j})^{\lbar (1-a_r)\lam_j}\right)
\]
The defining chain of $\upi$ is compatible with $\tw_0$, it follows (see Remark~\ref{modulodescription}) 
for the lifts $\tta_j$ of the $\tau_j$ that the $b_{\tta_j}^{\varpi_i}$ are generators of 
${}^t\ca_{\tw_0}^\xi\subset \cri_\xi[{}^tU^-]$. Note if $\tta\in W$ is a lift for $\tau\in W/W_\lam$, 
then $(b_\tau^{\varpi_j})^{\lam_j}=(b_\tta^{\varpi_j})^{\lam_j}$. 
So $b_\pi\in {}^t\ca_{\tw_0}^\xi$, and by reordering the elements we can write with the notation as 
in (\ref{generatorforadap}):
\[
m_\pi={\it c} b_1^{r_1}\cdots b_N^{r_N},
\]
where {\it c} is some root of unity. If all the exponents
$r_j$ are divisible by $\lbar$, then we know by Proposition~\ref{frobeniusminors} that
\[
p_\pi={\Fr'}^*(m_\pi) ={\it c}{\Fr'}^*(b_1^{r_1}\cdots b_N^{r_N})=
{\it c} b_1^{r_1/\lbar}\cdots b_N^{r_N/\lbar}
\]
is an element of $S_{\tw_0}$, the spanning set for $\ca_{\tw_0}^1$. 

To prove the divisibility property, let $\pi=(\tau_0,\ldots,\tau_r;a_1,\ldots,a_r)$
be an \hbox{L-S} path of shape $\lam$ compatible with $\tw_0$ and consider the following projection: 
write $\lam$ as a linear combination $\sum_\al \lam_\al \varpi_\al$ of fundamental
weights, fix a simple root $\al$, let $\tbar$ be the class of $\tau$ in $W/W_{\varpi_\al}$, and set 
\[
\pi_\al=(\tbar_0,\ldots,\tbar_r;a_1,\ldots,a_r)=(\sigma_0,\ldots,\sigma_t;c_1,\ldots,c_t).
\]
where the second pair of sequences is obtained from the first by omitting repetitions
among the $\tbar_j$ as well as the corresponding $b_j$, i.e., if 
\[
\sigma_{0}=\tbar_{0}\ldots=\tbar_{j-1}>\sigma_{1}=\tbar_j=\ldots=\tbar_{k-1}
>\sigma_{2}=\tbar_{k}=\ldots=\tbar_{m-1}>\ldots
\]
then $c_1=a_j, c_{2}=a_{k}, c_{3} = a_{m}, \ldots$. Set
$$
m_{\pi_\al}= \left(b_{\sigma_0}^{\varpi_\al}\right)^{\lbar\lam_\al c_1}
\left(b_{\sigma_1}^{\varpi_\al}\right)^{\lbar\lam_\al (c_2-c_1)}\cdots
\left(b_{\sigma_t}^{\varpi_\al}\right)^{\lbar\lam_\al (1-c_t)}\in \cri_\xi[{}^tU^-]
$$
Then $m_\pi=\xi^p m_{\pi_{\al_1}} m_{\pi_{\al_2}}\cdots m_{\pi_{\al_n}}\in \cri_\xi[{}^tU^-]$ for some
power $p$ of $\xi$. So the divisibility property, and hence the first part of the theorem, is a consequence 
of the following lemma, and the fact that $p_\pi$ and the canonical basis element $b_\pi$ coincide
(up to multiplication by an appropriate root of unity) follows from Proposition~\ref{lnt}.

\begin{lem} 
Set $c_0=0$ and $c_{r+1}=1$. The exponents have the property
$\lam_\al (c_{j+1}-c_{j})\in\bz$ for all $j=0,\ldots,r$.
\end{lem}
\par
\noindent
{\it Proof.}
Let $\pi_\al=(\tbar_0,\ldots,\tbar_r;a_1,\ldots,a_r)=(\sigma_0,\ldots,\sigma_t;c_1,\ldots,c_t)$ be as above
and suppose
\[
\sigma_{i-1}=\tbar_{j-1}>\sigma_{i}=\tbar_j=\tbar_{j+1}=\ldots=\tbar_{k-1}
>\sigma_{i+1}=\tbar_{k}
\]
Corresponding to the fixed reduced decomposition $\tw_0$ let 
$y_j=s_{i_1}\cdots s_{i_j}$ be as in (\ref{generatorforadap}). By assumption,
$\pi$ is compatible with $\tw_0$, so $\{\tau_0,\ldots,\tau_r\}\subset \{y_1,\ldots, y_N\}$.
It follows: there exist indices such that 
\[
\tau_k\le y_{t}<y_{t}s_\al\le\tau_{k-1}<\ldots<\tau_j\le y_{v}<y_{v}s_\al\le \tau_{j-1}
\]
and $y_{t}s_\al$ and $y_{v}s_\al$ are elements of $\{y_1,\ldots, y_N\}$. Set $\beta=y_t(\al)$,
so $s_\beta (y_t s_\al)=y_{t}$. The definition of an L-S path implies:
\[
c_{i+1}\lam_\al= a_k\lam_\al = a_k\langle \lam,\al^\vee \rangle
 = a_k\langle  y_{t}(\lam),\beta^\vee \rangle \in\bz.
\]
Similarly, set $\delta=y_v(\al)$,
so $s_\delta(y_v s_\al)=y_{v}$. The definition of an L-S path implies:
\[
c_{i}\lam_\al= a_j\lam_\al = a_j\langle \lam,\al^\vee \rangle 
 = a_j\langle  y_{v}(\lam),\delta^\vee \rangle \in\bz.
\]
It follows: $\lam_\al(c_{i+1}-c_{i})\in\bz$.
\endpf

%%%%%%%%%%%%%%%%%%%%%%%%%%%%%%%%%%%%%%%%%%%%%%%%%%%%%%%%%%%%%%%%%%%%%%%%%%%%%%%%%%%%%%%%%%%%%%%%%%%%%%%%%%%%
%%%%%%%%%%%%%%%%%%%%%%%%%%%%%%%%%%%%%%%%%%%%%%%%%%%%%%%%%%%%%%%%%%%%%%%%%%%%%%%%%%%%%%%%%%%%%%%%%%%%%%%%%%%%
\section{Examples}\label{examples}
%%%%%%%%%%%%%%%%%%%%%%%%%%%%%%%%%%%%%%%%%%%%%%%%%%%%%%%%%%%%%%%%%%%%%%%%%%%%%%%%%%%%%%%%%%%%%%%%%%%%%%%%%%%%
%%%%%%%%%%%%%%%%%%%%%%%%%%%%%%%%%%%%%%%%%%%%%%%%%%%%%%%%%%%%%%%%%%%%%%%%%%%%%%%%%%%%%%%%%%%%%%%%%%%%%%%%%%%%
In this section, we study some examples which illustrate our main theorem.
Let $\km$ be semisimple Lie algebra, let $\cb\subset U_q(\Ln^-)$ be the canonical
basis and denote $L_q(\Ln^-)$ the $\bc[q]$--submodule of $U_q(\Ln^-)$ spanned by the elements 
of the canonical basis $\cb$. Let $(B(\infty), L_q(\Ln^-))$ be the crystal of 
$U_q(\Ln^-)$ (see for example \cite{Ka}, or \cite{joseph}), consisting of 
the $\bc[q]$--module $L_q(\Ln^-)$, the basis
$$
B(\infty)=\{ b\bmod L_q(\Ln^-)\mid b\in\cb\},
$$
and the Kashiwara operators $\{f_\al,e_\al\mid\al\hbox{{\rm\ a\ simple\ root\/}}\}$ 
on $B(\infty)\cup\{0\}$.

In the same way as the dual canonical basis $\cb^*$
is compatible with the representations $V_q(\lam)^*$, the crystal
of $U_q(\Ln^-)$ can be seen as the limit of the crystals of the representations
$V_q(\lam)$. Recall (section~\ref{pathmodel}) that the path model can also be understood as a combinatorial model for the 
crystal of the representations. 

\begin{exam}\label{slthree}\rm
Let $\km$ be of type ${\tt A}_2$, denote $\al_1,\al_2$ the simple roots, then every
element $b\in\cb^*\subset\bc_1[U^-]$ of the dual canonical basis is the image
of some appropriately chosen standard monomial $p_\upi$.

\smallskip\noindent
{\it Proof.\/} Let $u_\infty\in B(\infty)$ be the highest weight element of $B(\infty)$ and let
$u_b$ be the class of $b$ in $B(\infty)$. 
The element $u_b\in B(\infty)$ can be written as
$$
u_b=f^\ell_{\al_1} f^m_{\al_2}f^n_{\al_1} u_\infty,\ \hbox{\rm\ where\ }m\ge n,
\ \hbox{\rm\ or\ }u_b=f^r_{\al_2} f^s_{\al_1}f^t_{\al_2} u_\infty\ \hbox{\rm\ where\ }s\ge t,
$$
and $s=\ell+n$, $r=\max\{n,m-\ell\}$ and $t=\min\{m-n,\ell\}$ (see \cite{Li$_4$} or \cite{berzev}). 

Note if $u_b=f^\ell_{\al_1} f^m_{\al_2}f^n_{\al_1} u_\infty$ is such that $\ell<m-n$, then in 
the corresponding expression $u_b=f^r_{\al_2} f^s_{\al_1}f^t_{\al_2} u_\infty$ we have 
$r=m-\ell>n= \ell+n - \ell = s-t$. So in the following we just write $u_b=f^\ell_{\al_i} f^m_{\al_j}f^n_{\al_i} u_\infty$
and suppose that $i$ and $j$ are chosen such that $\ell\ge m-n\ge 0$. 

Set $\lam=(\ell+2n-m)\varpi_i+(m-n)\varpi_j$, and let $\upi_\ulam=(id,id,id)$ be the unique
standard tuple of shape $\ulam=(n\varpi_i\,,(m-n)\varpi_j\,,(\ell+n-m)\varpi_i\,)$ of weight $\upi_\ulam(1)=\lam$. 
By applying the root operators (see \cite{Li$_3$}), we get for $\upi=f^\ell_{\al_i} f^m_{\al_j}f^n_{\al_i} \pi_\ulam$ the 
following standard sequence of shape $\ulam$:
$$
\upi=\big(\pi_1=(s_j s_i)\,,\pi_2=(s_i s_j)\,,\pi_3=(s_i)\big).
$$
Note that the sequence $(s_i s_j s_i\,, s_i s_j\,, s_i)$ is a defining chain for $\upi$
which is obviously compatible with the reduced decomposition $w_0=s_i s_j s_i$. It follows that
the image of the standard monomial $p_\upi=p_{s_j s_i} p_{s_j s_i} p_{s_i}$ of shape $\ulam$
in $\bc[U^-]$ is, up to multiplication by a root of unity, the element $b=b_\upi$ of the 
dual canonical basis $\cb^*$.
\end{exam}

\begin{exam}\label{slfour}\rm
The case $\km={\mathfrak sl}_3$ is very special, in general it is not possible to find for every element in $B(\infty)$ an ``appropriate'' standard tuple
of L--S paths with a defining chain compatible with a reduced decomposition of $w_0$.
An example is for $\km={\mathfrak sl}_4$ the element $b\in\cb^*$ such that for its class
in $B(\infty)$ holds $u_b=f_{\al_2} f_{\al_1} f_{\al_3}u_\infty$.
Nevertheless, it is easy to see by Theorem~\ref{pathandcanonicalzigzag}, that
the other elements of $B(\infty)$ of the same weight
$$
u_{b_1}=f_{\al_1} f_{\al_2} f_{\al_3}u_\infty\,,\quad u_{b_2}=f_{\al_3} f_{\al_2} f_{\al_1}u_\infty\,,\quad
u_{b_3}=f_{\al_3} f_{\al_1} f_{\al_1}u_\infty
$$
correspond all to elements of the dual canonical basis which are ``standard monomials'':
$b_1=p_{(s_1s_2s_3)}\in H^0(G/B,\cl_{\varpi_3})$, $b_2=p_{(s_3s_2s_1)}\in H^0(G/B,\cl_{\varpi_1})$
and $b_3=p_{(s_1s_3s_2)}\in H^0(G/B,\cl_{\varpi_2})$.

But let now $\lam=\varpi_1+\varpi_3$, then there exists an L--S path of shape $\lam$ that corresponds 
under the identification between crystal graphs and path models to $u_b$, it 
is the path $\pi=(s_2 s_3 s_1, s_3s_1; \frac{1}{2})$. The $\km$-module $V(\lam)$ is isomorphic to $\km$ endowed with the adjoint action. Via this isomorphism, the canonical basis of $V(\lam)$ is given (up to scalar multiples) by a Chevalley basis. As a function, $p_\pi$
is the dual of the canonical basis element $h_{\al_2}$, thus $b=p_\pi$ can again be constructed
as a ``standard monomial'' for an appropriate choice of $\lam$, but this time not in the context 
of Theorem~\ref{pathandcanonicalzigzag}. 
\end{exam}

\begin{exam}\rm We suppose again $\km={\mathfrak sl}_3$. Set $\lambda=2\varpi_1+2\varpi_2$. We shall
explicit the $p_\pi$ corresponding to the zero weight space $V(\lambda)_0$
of  $V(\lambda)$.\par
Let's introduce some notation. Let
$$a=p_{s_1}^{\varpi_1},\,b=p_{s_2}^{\varpi_2},\,c=p_{s_2s_1}^{\varpi_1},\,d=
p_{s_1s_2}^{\varpi_2},$$
be the elements of $\bc[U^-]$ as defined in section~\ref{pathmodel}. It is well
known (see \cite{Cal}), that the dual canonical basis is given by
$\{a^kb^lc^sd^t,\,k,l,s,t\in\bn,\, kl=0\}$. Moreover, it is easy to
see that $ab=c+d$. Indeed, $ab$ decomposes in the base $\{c,d\}$ of the
subspace of  $\bc[U^-]$ of weight $\alpha_1+\alpha_2$. Both
coefficients in this base are equal by symmetry and equal to one by
\cite{Cal2}, Theorem 2.3.
\par
The LS paths corresponding to $V(\lambda)_0$ are
$$
\pi_1=(s_2s_1,s_1;\, {1\over 2}\,),\,\pi_2=(w_0,\hbox{Id};\,
{1\over 2}\,),\,\pi_3=(s_1s_2,s_2;\,{1\over 2}\,).
$$
Remark that $\pi_2$ is compatible with a reduced decomposition while $\pi_1$
and $\pi_3$ are not. As in the proof of Theorem 9, with $l=2$, we obtain, up
to a root of one,
$$p_{\pi_1}=cba,\;p_{\pi_2}=cd,\;p_{\pi_3}=dab.$$
Hence,
$$p_{\pi_1}=c^2+cd,\;p_{\pi_2}=cd,\;p_{\pi_3}=cd+d^2.$$
Remark that $p_{\pi_2}$ is an element of the dual canonical basis while
$p_{\pi_1}$ and $p_{\pi_3}$ are not. Under the identification of the path model
with crystal basis elements we have in the dual canonical basis:
$b_{\pi_1}=c^2$, $b_{\pi_2}=cd $ and  $b_{\pi_3}= d^2$, so the transformation
matrix reads as: $p_{\pi_1}= b_{\pi_1} + b_{\pi_2}$, $p_{\pi_2}= b_{\pi_2}$
and $p_{\pi_3}= b_{\pi_3} + b_{\pi_2}$.
\end{exam}

\begin{exam}\rm We suppose again $\km={\mathfrak sl}_3$, set now
$\lambda=m\varpi_1+m\varpi_2$. To simplify the notation of the L--S paths, 
we still write $\pi=(\ut,\ua)$ for $\ua=(a_1,a_2,a_3)$, even if we
have equality in $0\le a_1\le a_2\le a_3\le 1$. In case of 
equality, the reader has to omit the corresponding parameters and Weyl group elements
to avoid double counting. 

We claim that the $p_\pi$ have a simple expression
in terms of the canonical basis: if $\pi=(\ut,\ua)$ is compatible
with some reduced decomposition, then $p_\pi=b_{\pi}$ is an
element of the dual canonical basis (after a possible renormalization
with a root of unity). 

If $\pi=(\ut,\ua)$ is not compatible 
with a reduced decomposition, then let $\ell$ be even and minimal
with the property that $\ell a_i\in\bn$ for all $i$.
Again, after a possible renormalization (multiplication
with a root of unity):

\begin{prop}
If $\pi=(\ut,\ua)$ is not compatible 
with a reduced decomposition, then set $t=\min\{m(a_2-a_1),m(a_3-a_2)\}$. The following
transition relation holds between the canonical basis and the path vectors:
$$
p_{(\ut;\ua)}=
\sum_{j=0}^{\lfloor t\rfloor}{t\ell \choose j\ell}  b_{(\ut;a_1+{\scriptscriptstyle\frac{j}{m}},a_2,a_3-\frac{j}{m})}
$$
\end{prop}
\vskip 8pt
{\it Proof.\/}
To prove the statement, we divide the set of L--S paths in groups. 
The first group is the set $B(\lam)^1=B(\lam)^{1,1}\cup B(\lam)^{1,2}$, where:
\[
\begin{array}{c}
B(\lam)^{1,1}=\{\pi=(w_0,s_1s_2,s_1,id;a_1,a_2,a_3 )\mid m a_1,m a_2, ma_3\in\bn\}\cr
B(\lam)^{1,2}=\{\pi=(w_0,s_2s_1,s_2,id;a_1,a_2,a_3 )\mid m a_1,m a_2, ma_3\in\bn\}\cr
\end{array}
\]
These paths are all compatible with a reduced decomposition, so the elements $p_\pi$
are canonical basis elements: Set $x_1=a_1$, $x_2=a_2-a_1$, $x_3=a_3-a_2$, then, after renormalizing
the $p_\pi$ (multiplication by a root of unity), we have (with same notation as above):
\[
\begin{array}{c}
b_\pi= p_\pi= a^{m(x_2+x_3)} c^{m x_1} d^{m(x_1+x_2)}\hbox{\rm\ for\ }\pi\in B(\lam)^{1,1}\cr
b_\pi= p_\pi= b^{m(x_2+x_3)} c^{m (x_1+x_2)} d^{mx_1}\hbox{\rm\ for\ }\pi\in B(\lam)^{1,2}
\end{array}
\]
In other words, let $a^uc^vd^w$ (respectively $b^u c^v d^w$) be an element of the
dual canonical basis such that 
$$
v\le w\le u+v\le m\quad (\hbox{{\rm respectively:\ }}w\le v\le u+w\le m)
$$
In both cases, this canonical basis element is equal to $p_\pi$, where 
$$
\pi=\left\{
\begin{array}{ll}
(w_0,s_1s_2,s_1,id;\frac{v}{m},\frac{w}{m},\frac{u+v}{m})&\hbox{{\rm\ for\ }}a^uc^vd^w\cr
(w_0,s_2s_1,s_2,id;\frac{w}{m},\frac{v}{m},\frac{u+w}{m})&\hbox{{\rm\ for\ }}b^u c^v d^w\hbox{{\rm\ respectively}}\cr
\end{array}
\right.
$$

The second set of L--S paths of shape $\lam$ is in this example the union 
$B(\lam)^2=B(\lam)^{2,1}\cup B(\lam)^{2,2}$, where:
\[
\begin{array}{c}
B(\lam)^{2,1}=\{\pi=(w_0,s_2s_1,s_1,id;a_1,a_2,a_3 )\mid m a_1,2m a_2, ma_3\in\bn\}\cr
B(\lam)^{2,2}=\{\pi=(w_0,s_1s_2,s_2,id;a_1,a_2,a_3 )\mid m a_1,2m a_2, ma_3\in\bn\}\cr
\end{array}
\]
By construction, we get (up to multiplication by a root of unity): Let $\ell$ be minimal
and even such that $\ell a_i\in\bn$ for all $i$, then, for $\pi\in B(\lam)^{2,1}$, we have up
to multiplication by a root of unity $\tilde{r}$:
$$
p_\pi=
\left\{
\begin{array}{rl}
\tilde{r}{\Fr'}^*(c^{mx_1\ell}d^{m x_1\ell}c^{m x_2\ell}(ab)^{m x_3\ell}b^{m(x_2-x_3)\ell})&\hbox{\rm\ if\ }x_2\ge x_3\cr
\tilde{r}{\Fr'}^*(c^{mx_1\ell}d^{m x_1\ell}c^{m x_2\ell}(ab)^{m x_2\ell}a^{m(x_3-x_2)\ell})&\hbox{\rm\ if\ }x_3\ge x_2
\end{array}
\right.
$$
and hence by Proposition~\ref{frobeniusminors} for some roots of unity $r'$:
$$
p_\pi=
\left\{
\begin{array}{rl}
r'b^{m(x_2-x_3)} c^{mx_1}d^{m x_1}{\Fr'}^*\big(c^{m x_2\ell}(ab)^{m x_3\ell}\big)&\hbox{\rm\ if\ }x_2\ge x_3\cr
r'a^{m(x_3-x_2)} c^{mx_1}d^{m x_1}{\Fr'}^*\big(c^{m x_2\ell}(ab)^{m x_2\ell}\big)&\hbox{\rm\ if\ }x_3\ge x_2.
\end{array}
\right.
$$
Replacing $ab$ by $(c+d)$ (note that $c$ and $d$ commute) we get after renormalizing
(by a root of unity):
$$
p_\pi=
\left\{
\begin{array}{ll}
\sum_{j=0}^{mx_3\ell}{m\ell x_3 \choose j} b^{m(x_2-x_3)} c^{mx_1}d^{m x_1}{\Fr'}^*(c^{m (x_2+x_3)\ell-j}d^{j})&\hbox{\rm\ if\ }x_2\ge x_3\cr
\sum_{j=0}^{mx_2\ell}{m\ell x_2 \choose j} a^{m(x_3-x_2)} c^{mx_1}d^{m x_1}{\Fr'}^*(c^{m 2x_2\ell-j}d^j)&\hbox{\rm\ if\
}x_3\ge x_2
\end{array}
\right.
$$
Since $2mx_2$ respectively $m(x_2+x_3)$ are integers, it follows that ${\Fr'}^*$ maps a summand to zero
for $j$ not divisible by $\ell$, and hence:
$$
p_\pi=
\left\{
\begin{array}{ll}
\sum_{j=0}^{\lfloor mx_3\rfloor}{m\ell x_3 \choose j\ell} b^{m(x_2-x_3)} c^{mx_1}d^{m x_1}c^{m (x_2+x_3)-j}d^{j}&\hbox{\rm\ if\ }x_2\ge x_3\cr
\sum_{j=0}^{\lfloor mx_2\rfloor}{m\ell x_2 \choose j\ell} a^{m(x_3-x_2)} c^{mx_1}d^{m x_1}c^{m 2x_2-j}d^j&\hbox{\rm\ if\
}x_3\ge x_2
\end{array}
\right.
$$
and thus, after reordering, we get for $\pi\in B(\lam)^{2,1}$:
\begin{eqnarray}\label{firstsum}
p_\pi=
\left\{
\begin{array}{ll}
\sum_{j=0}^{\lfloor mx_3\rfloor}{m\ell x_3 \choose j\ell}  b^{m(x_2-x_3)} c^{m(x_1+x_2+x_3)-j}d^{m x_1+j}&\hbox{\rm\ if\ }x_2\ge x_3\cr
\sum_{j=0}^{\lfloor mx_2\rfloor}{m\ell x_2 \choose j\ell}  a^{m(x_3-x_2)} c^{m(x_1+2x_2)-j}d^{m x_1+j}&\hbox{\rm\ if\ }x_3\ge x_2
\end{array}
\right.
\end{eqnarray}
One sees easily that in terms of the root operators one has 
$$
\pi=f_1^{m(x_1+2x_2)} f_2^{m(2x_1+x_2+x_3)} f_1^{mx_1}(id),
$$ 
so, by Example~\ref{slthree}, the corresponding element $b_\pi$ in the dual canonical basis is 
$$
b_\pi=
\left\{
\begin{array}{ll}
 b^{m(x_2-x_3)} c^{m(x_1+x_2+x_3)}d^{m x_1}&\hbox{\rm\ if\ }x_2\ge x_3\cr
 a^{m(x_3-x_2)} c^{m(x_1+2x_2)}d^{m x_1}&\hbox{\rm\ if\ }x_3\ge x_2.
\end{array}
\right.
$$
i.e., $b_\pi$ is the summand for $j=0$ in the expression above.
Note that the exponents in the summands for $p_\pi$ in (\ref{firstsum}) satisfy all the following inequalities:
$$
\begin{array}{cccccccl}
 2w&\le& u+v+w&\le& 2v   &\le& 2m&\hbox{\rm\ for\ } b^uc^vd^w\cr
 2w&\le& v+w  &\le& 2(u+v)&\le& 2m&\hbox{\rm\ for\ } a^uc^vd^w.
\end{array}
$$
On the other hand, any triple $(u,v,w)$ of integers satisfying these inequalities gives an element of
the dual canonical basis of $V(\lam)^*$ and is hence equal to $b_\eta$ for some L--S path $\eta$ of shape
$\lam$. We obtain the corresponding path by the following rules:
$$
\eta=\left\{
\begin{array}{ll}
(w_0,s_2s_1,s_1,id;\frac{w}{m},\frac{u+v+w}{2m},\frac{v}{m})&\hbox{{\rm\ for\ }}b_\eta=b^u c^v d^w\cr
(w_0,s_2s_1,s_1,id;\frac{w}{m},\frac{v+w}{2m},\frac{u+v}{m})&\hbox{{\rm\ for\ }}b_\eta=a^u c^v d^w\cr
\end{array}
\right.
$$
In the same way one gets for $\pi\in B(\lam)^{2,2}$:
\begin{eqnarray}\label{secondsum}
p_\pi=
\left\{
\begin{array}{ll}
\sum_{j=0}^{\lfloor mx_3\rfloor}{m\ell x_3 \choose j\ell}  a^{m(x_2-x_3)} c^{m x_1+j}d^{m(x_1+x_2+x_3)-j}&\hbox{\rm\ if\ }x_2\ge x_3\cr
\sum_{j=0}^{\lfloor mx_2\rfloor}{m\ell x_2 \choose j\ell}  b^{m(x_3-x_2)} c^{m x_1+j}d^{m(x_1+2x_2)-j}&\hbox{\rm\ if\
}x_3\ge x_2
\end{array}
\right.
\end{eqnarray}
As above, in terms of the root operators one has 
$$
\pi=f_2^{m(x_1+2x_2)} f_1^{m(2x_1+x_2+x_3)} f_2^{mx_1}(id),
$$ 
so the corresponding element $b_\pi$ in the dual canonical basis is the summand for $j=0$ 
in the expression above:
$$
b_\pi=
\left\{
\begin{array}{ll}
 a^{m(x_2-x_3)} c^{m x_1} d^{m(x_1+x_2+x_3)}&\hbox{\rm\ if\ }x_2\ge x_3\cr
 b^{m(x_3-x_2)} c^{m x_1} d^{m(x_1+2x_2)}&\hbox{\rm\ if\ }x_3\ge x_2
\end{array}
\right.
$$
The exponents in the expression for $p_\pi$ in (\ref{secondsum}) satisfy the following inequalities:
$$
\begin{array}{cccccccl}
 2v&\le& u+v+w&\le &2w    &\le &2m&\hbox{\rm\ for\ } b^uc^vd^w\cr
 2v&\le& v+w  &\le &2(u+w)&\le &2m&\hbox{\rm\ for\ } a^uc^vd^w
\end{array},
$$
and for any triple $(u,v,w)$ of integers satisfying these inequalities, we can obtain the corresponding 
path by the following rules:
$$
\eta=\left\{
\begin{array}{ll}
(w_0,s_1s_2,s_2,id;\frac{v}{m},\frac{u+v+w}{2m},\frac{w}{m})&\hbox{{\rm\ for\ }}b_\eta=a^u c^v d^w\cr
(w_0,s_1s_2,s_2,id;\frac{v}{m},\frac{v+w}{2m},\frac{u+v}{m})&\hbox{{\rm\ for\ }}b_\eta=b^u c^v d^w\cr
\end{array}
\right.
$$
As a consequence we can now describe in terms of the path model the sum expressions in (\ref{firstsum})
and (\ref{secondsum}): For all $\pi=(\ut,\ua)\in B(\lam)^2$, $\ua=(a_1,a_2,a_3)$, we get for
$t=\min\{mx_2,mx_3\}$
\begin{eqnarray}\label{firstfinsum}
p_{(\ut;\ua)}=
\sum_{j=0}^{\lfloor t\rfloor}{t\ell \choose j\ell}  b_{(\ut;a_1+{\scriptscriptstyle\frac{j}{m}},a_2,a_3-\frac{j}{m})}
\end{eqnarray}
\end{exam}
\endpf

\goodbreak
\parindent 0pt
P.C.: Institut Girard Desargues, UPRES-A-5028

Universit\'e Claude Bernard, Lyon I, Bat 101

69622 Villeurbane Cedex, France

email: caldero@desargues.univ-lyon1.fr

\smallskip
P.L: Fachbereich Mathematik 

Universit\"at Wuppertal, Gau\ss-Stra\ss e 20 

42097 Wuppertal, Germany

email: littelmann@math.uni-wuppertal.de

\smallskip
Keywords: Quantum groups, Canonical basis, Standard Monomial Theory
%%%%%%%%%%%%%%%%%%%%%%%%%%%%%%%%%%%%%%%%%%%%%%%%%%%%%%%%%%%%%%%%%%%%%%%%%%%%%%%%%%%%%
%%%%%%%%%%%%%%%%%%%%%%%%%%%%%%%%%%%%%%%%%%%%%%%%%%%%%%%%%%%%%%%%%%%%%%%%%%%%%%%%%%%%%
%%%%%%%%%%%%%%%%%%%%%%%%%%%%%%%%%%%%%%%%%%%%%%%%%%%%%%%%%%%%%%%%%%%%%%%%%%%%%%%%%%%%%
\end{document}